\documentclass{article}
\usepackage{amsthm}
\usepackage{amsmath}
\usepackage{amssymb}
\usepackage{mathrsfs}
\usepackage{tikz}
\usepackage{graphicx}
\usepackage{relsize}
\usepackage{geometry}
\usepackage{enumitem}
\usepackage{xcolor}
\usepackage{booktabs}
\usepackage{authblk}
\usepackage{float}
\usepackage[hidelinks]{hyperref}

\usetikzlibrary{arrows.meta}
\newtheorem{theorem}{Theorem}
\newtheorem{proposition}{Proposition}
\newtheorem{lemma}{Lemma}
\newtheorem{corollary}{Corollary}
\theoremstyle{definition}
\newtheorem{definition}{Definition}

\newtheorem{remark}{Remark}
\newtheorem{example}{Example}

\newcommand{\N}{\mathbb N_0}
\newcommand{\im}{\operatorname{im}}
\newcommand{\rk}{\operatorname{rank}}
\newcommand{\id}{\operatorname{id}}

\title{Division rings of fractions of LV-algebras}

\author[1]{Huynh Viet Khanh}
\author[1]{Vo Hoang Thanh}
\affil[1]{Department of Mathematics and Informatics, HCMC University of Education, Ho Chi Minh City, Vietnam\\\texttt{khanhhv@hcmue.edu.vn}, \texttt{thanhvh@hcmue.edu.vn}}
\date{ }

\begin{document}

\maketitle

\begin{abstract}
Let $K$ be a field and let $A=L_K(E,\omega)$ be an LV-algebra of a row-finite weighted graph. We prove that $A$ embeds in a division ring if and only if $E$ has one vertex and its ordered loop weights satisfy $w_n=n$ and $w_j>j$ for $1\le j<n$. Under these conditions, $A$ has a universal division ring of fractions, which we identify with the free field on $2\sum_{j=1}^n w_j-n^2$ generators over $K$. We construct a free generating set explicitly by factoring the matrix of generators. We classify the resulting division rings over $K$ and bound the stabilization required to recognize invertible matrices. These results include mixed weights, for which the LV-algebra need not be a semifir.
\end{abstract}

{\small

\par}

\section{Introduction}
\label{sec:introduction}

The absence of zero-divisors is necessary for a ring to embed in a division ring. As Cohn notes in \cite[\S1.2, p.~9]{C}, a noncommutative domain need not admit such an embedding. We study this distinction for LV-algebras, the weighted Leavitt path algebras of graphs satisfying the LV conditions of Hazrat and Preusser \cite[Definition~38, p.~1077]{HP}. Throughout, $K$ denotes a commutative field. Let $A=L_K(E,\omega)$ be an LV-algebra with $E$ row-finite. The vertices of $E$ are nonzero in $A$ by Preusser's normal form \cite[Remark~7(e), p.~129]{P}. Distinct vertices are orthogonal, so $A$ can be a domain only if $E$ has one vertex. In that case $E$ is a finite LV-rose by row-finiteness, and $A$ is a domain by Hazrat and Preusser \cite[Theorem~41, p.~1079]{HP}. Example~\ref{ex:domain-no-embedding} exhibits a domain with weights $(2,2,4,4)$ that does not embed in a division ring. We determine which weights permit an embedding and describe the division ring of fractions distinguished by a universal property.

For a ring $R$ and a set $\Sigma$ of matrices over $R$, Cohn constructs the universal localization $R_\Sigma$ by adjoining inverses to these matrices \cite[Theorem~4.1.3, p.~156]{C}. The localization is universal for homomorphisms which invert these matrices. The essential step is to choose $\Sigma$ so that $R_\Sigma$ is a division ring and $R\to R_\Sigma$ is injective. For a right hereditary ring with a faithful integer-valued projective rank function, it follows from Cohn's results on rank and localization \cite[Propositions~4.6.11 and~4.6.13, pp.~193--195]{C} that one may take $\Sigma$ to consist of the square matrices which are full with respect to that rank. The localization preserves the resulting matrix ranks. Thus information about finitely generated projective modules can establish both the existence of an embedding and an intrinsic construction of a containing division ring.

For a finite weighted rose, Preusser shows that $A$ is hereditary and describes $V(A)$ in \cite[Theorem~14, p.~134]{P}. For such a rose, every embedding $\lambda\colon A\hookrightarrow F$ induces a faithful integer-valued projective rank function defined by $\rho_F([P])=\dim_F(P\otimes_AF)$. The relations in $V(A)$ therefore impose numerical conditions on the weights. Conversely, when these conditions define a faithful rank function, Cohn's localization constructs an embedding. More precisely, Theorem~\ref{thm:embedding} states that $A$ embeds in a division ring if and only if $E$ has one vertex and the weights of its $n$ structured loops, ordered as $w_1\le\cdots\le w_n$, satisfy
$$
w_n=n\quad\text{and}\quad w_j>j\quad\text{for }1\le j<n.
$$
Under these conditions, the projective rank function $\rho$ on $A$ is unique, integer-valued and faithful.

Mixed weights explain the role of projective rank functions in the construction. Corollary~\ref{cor:semifir} shows that an LV-algebra is a semifir precisely when its graph is a rose with $n$ loops, all of weight $n$. The rose with weights $(2,3,3)$ satisfies the embedding criterion but does not define a semifir. Example~\ref{ex:mixed} also exhibits a full square matrix over this algebra which is not stably full. Thus the matrices selected by $\rho$ cannot in general be described by fullness alone. In the semifir case, Cohn's construction \cite[Corollary~4.5.9, p.~182]{C} inverts every full square matrix.

A further question is which division ring of fractions to choose. As Cohn observes in \cite[\S4.1, p.~156]{C}, even division rings generated by the image of a given ring need not be isomorphic. In Cohn's terminology \cite[\S4.4, p.~176]{C}, a universal $A$-field specializes to every division ring generated by a homomorphic image of $A$. Let $\Sigma_\rho$ be the set of square matrices over $A$ which are full with respect to $\rho$. In Theorem~\ref{thm:universal}, we prove that the localization $D_K(E,\omega)=A_{\Sigma_\rho}$ is a universal division ring of fractions of $A$, unique up to a unique $A$-isomorphism. The proof uses the uniqueness of $\rho$ to compare singular kernels and then applies Cohn's specialization theorem \cite[Theorem~4.4.1, p.~172]{C}.

For $r\ge1$, let $\mathcal F_K(r)$ denote the universal division ring of fractions of the free associative $K$-algebra on $r$ generators. Theorem~\ref{thm:free-coordinates} proves
$$
D_K(E,\omega)\cong_K\mathcal F_K(r),\quad\text{where}\quad r=2\sum_{j=1}^n w_j-n^2.
$$
Lemma~\ref{lem:gauss-factors} constructs the factorization $X=LPT^{-1}$ of the matrix of generators from Section~\ref{sec:preliminaries}. The elements $p_j$ for $1\le j\le n$ and $l_{ij},t_{ji}$ for $1\le j<i\le w_j$ form a free generating set of $D_K(E,\omega)$ over $K$. Their number is $n+2\sum_j(w_j-j)=r$. We construct a homomorphism from $A$ to a free field on corresponding generators and use Cohn's specialization theorem to obtain the required isomorphism. The identification therefore uses the relations among the graph generators as well as the computation of $V(A)$.

According to Cohn \cite[Corollary~5.8.14, p.~258]{C}, the cardinality of a free generating set is an invariant of a free field over $K$. In view of the isomorphism above, $r$ classifies the resulting division rings up to $K$-isomorphism. Proposition~\ref{prop:attainable-ranks} determines all values of $r$ which occur. Example~\ref{ex:same-division-different-algebras} exhibits nonisomorphic LV-algebras with $K$-isomorphic universal division rings. In Example~\ref{ex:same-monoid-different-division}, there is a monoid isomorphism $\varphi\colon V(A)\to V(A')$ satisfying $\varphi([A])=[A']$, although the universal division rings are not $K$-isomorphic. Thus the projective monoid detects the rank obstruction to embedding, while the explicit free generators determine the resulting division ring.

The matrix localization also permits a description of arithmetic in $D_K(E,\omega)$. Proposition~\ref{prop:stablefull} characterizes invertibility over this division ring by stable fullness over $A$. Proposition~\ref{prop:bounded-stabilization} bounds the size of the identity block needed to test stable fullness, in terms of the weights. Proposition~\ref{prop:fractions} writes every element as $uM^{-1}v$ with $M\in\Sigma_\rho$, as in Cohn's matrix representation of fractions \cite[Theorem~4.2.1(d), p.~158]{C}. When $K$ is given by a presentation with solvable word problem, the explicit free generators allow us to use Cohn's result on the word problem for free fields \cite[Corollary~6.6.8, p.~323]{C} to decide equality of rational expressions and invertibility of square matrices over $A$ in $D_K(E,\omega)$.

Section~\ref{sec:preliminaries} introduces LV-algebras and recalls Preusser's description of their projective modules. Section~\ref{sec:criterion} proves the embedding criterion and characterizes the semifir case. Section~\ref{sec:construction} constructs the universal division ring and proves its specialization property. Section~\ref{sec:free-coordinates} determines free generators and classifies the resulting division rings over $K$. Section~\ref{sec:matrices} treats matrices and rational expressions.

\section{LV-algebras and their projective modules}
\label{sec:preliminaries}

A \emph{weighted graph} $(E,\omega)$ has vertex set $E^0$, structured edge set $E^{\mathrm{st}}$, source and range maps $s,r\colon E^{\mathrm{st}}\to E^0$, and a weight function $\omega\colon E^{\mathrm{st}}\to\mathbb N$, as in Preusser \cite[Definitions~1--2 and Remark~3, p.~127]{P}. We assume that $E$ is \emph{row-finite}, so $s^{-1}(v)$ is finite for each $v\in E^0$. A vertex $v$ is a \emph{sink} if $s^{-1}(v)=\varnothing$. The \emph{weight} of a vertex $v$ that is not a sink is defined by $\omega(v)=\max\{\omega(e):s(e)=v\}$. We use the LV conditions of Hazrat and Preusser \cite[Definition~38, p.~1077]{HP} and call $(E,\omega)$ an \emph{LV-graph} if it is connected, $E^{\mathrm{st}}\ne\varnothing$, $\omega(e)\ge2$ for every $e\in E^{\mathrm{st}}$, and
$$
\bigl|\{e\in s^{-1}(v):\omega(e)=\omega(v)\}\bigr|\ge2
$$
for every vertex $v$ that is not a sink. For each structured edge $e$, denote the corresponding edges by $e_1,\ldots,e_{\omega(e)}$ and their \emph{ghost edges} by $e_1^*,\ldots,e_{\omega(e)}^*$.

\begin{definition}
\label{def:wlpa}
As in Preusser \cite[Definition~4, p.~128]{P}, the \emph{weighted Leavitt path algebra} $L_K(E,\omega)$ is the universal $K$-algebra generated by the vertices together with the symbols $e_i,e_i^*$, where $e\in E^{\mathrm{st}}$ and $1\le i\le\omega(e)$, subject to the following relations.
\begin{enumerate}[label=\textup{(W\arabic*)}]
\item $uv=\delta_{uv}u$ for all $u,v\in E^0$.
\item $s(e)e_i=e_i=e_ir(e)$ and $r(e)e_i^*=e_i^*=e_i^*s(e)$ for every $e\in E^{\mathrm{st}}$ and $1\le i\le\omega(e)$.
\item For every vertex $v$ that is not a sink and all $1\le i,j\le\omega(v)$,
$$
\sum_{e\in s^{-1}(v)}e_ie_j^*=\delta_{ij}v.
$$
\item For every vertex $v$ that is not a sink and all $e,f\in s^{-1}(v)$,
$$
\sum_{i=1}^{\omega(v)}e_i^*f_i=\delta_{ef}r(e).
$$
\end{enumerate}
Here $e_i=e_i^*=0$ whenever $i>\omega(e)$. An \emph{LV-algebra} is $L_K(E,\omega)$ for an LV-graph $(E,\omega)$.
\end{definition}

The algebra $L_K(E,\omega)$ is unital if and only if $E^0$ is finite, in which case $1=\sum_{v\in E^0}v$.

A weighted graph with a single vertex is a \emph{weighted rose}. For the rest of this section, fix an LV-rose $(E,\omega)$ with vertex $v$ and put $A=L_K(E,\omega)$. Its structured edges are loops, and their number $n$ is finite with $n\ge2$. Order the loops $e^{(1)},\ldots,e^{(n)}$ so that their weights satisfy $w_1\le\cdots\le w_n$, and put $d=w_n$. Write $x_{ij}=e^{(j)}_i$ and $y_{ji}=(e^{(j)}_i)^*$ for $1\le i\le w_j$. Set $x_{ij}=y_{ji}=0$ for $i>w_j$, and form
$$
X=(x_{ij})\in M_{d\times n}(A)\quad\text{and}\quad
Y=(y_{ji})\in M_{n\times d}(A).
$$
Here $v=1_A$, and $A$ has the presentation
$$
A=K\langle x_{ij},y_{ji}:1\le j\le n,\ 1\le i\le w_j\rangle
\big/(XY-I_d,\ YX-I_n).
$$

All modules below are right modules. Vectors are columns, and matrices act on the left. Let $R$ be a unital ring. For a finitely generated projective $R$-module $P$, write $[P]$ for its isomorphism class, so that $[P]=[Q]$ if and only if $P\cong_R Q$. Addition on the set $V(R)$ of these classes is defined by
$$
[P]+[Q]=[P\oplus Q].
$$
Direct sum respects isomorphisms and is associative and commutative up to isomorphism. Thus this operation makes $V(R)$ a commutative monoid with identity $[0]$, the class of the zero module. In particular, $[R^n]=n[R]$ for every $n\in\N$.

The \emph{Grothendieck group} $K_0(R)$ is the \emph{group completion} of $V(R)$, with canonical map $\iota\colon V(R)\to K_0(R)$ \cite[\S4.6, p.~186]{C}. Every additive map $f\colon V(R)\to H$ to an abelian group factors uniquely as $f=\bar f\circ\iota$, where $\bar f\colon K_0(R)\to H$ is a group homomorphism.

We use Cohn's definition of a \emph{projective rank function} \cite[\S4.6, p.~187]{C}. Such a function on $R$ is a map $\rho\colon V(R)\to\mathbb R_{\ge0}$ satisfying
$$
\rho([0])=0,\qquad
\rho([P\oplus Q])=\rho([P])+\rho([Q])
\quad\text{and}\quad \rho([R])=1.
$$
Since $\rho$ is defined on isomorphism classes, $P\cong_R Q$ implies $\rho([P])=\rho([Q])$. We write $\rho(P)=\rho([P])$; in particular, $\rho(R^n)=n$. The rank is \emph{faithful} if $\rho(P)>0$ for every nonzero finitely generated projective $P$, and it is \emph{integer-valued} if $\rho(P)\in\N$ for every finitely generated projective $P$.

For the rose above, let $d_1<\cdots<d_k=d$ be the distinct weights, let $a_i$ be the number of loops of weight $d_i$, and put $N_i=a_1+\cdots+a_i$, with $d_0=N_0=0$ and $N_k=n$. The LV condition implies $d_1\ge2$ and $a_k\ge2$.

\begin{proposition}
\label{prop:monoid}
For the LV-rose $(E,\omega)$ fixed above, the algebra $A=L_K(E,\omega)$ is left and right hereditary, and there is an isomorphism
$$
V(A)\cong\bigl\langle v,q_1,\ldots,q_{k-1}\ \bigm|\
q_{i-1}+(d_i-d_{i-1})v=q_i+a_iv\quad\text{for }1\le i\le k\bigr\rangle,
$$
where $v=[A]$ and $q_0=q_k=0$. If $k>1$, then $q_i\ne0$ and $q_i\ne mv$ for every $1\le i<k$ and every $m\in\N$.
\end{proposition}

\begin{proof}
Since the rose is finite, it follows from Preusser's theorem \cite[Definition~11, p.~132, and Theorem~14, p.~134]{P} that $A$ is left and right hereditary and that $V(A)$ has the displayed presentation. To examine the generators, let $F$ be the free commutative monoid on the symbols $v,q_1,\ldots,q_{k-1}$. Let $\ell\colon F\to\N$ be the function defined by
$$
\ell\bigl(cv+\sum_{j=1}^{k-1}b_jq_j\bigr)=c+\sum_{j=1}^{k-1}b_j,
\qquad\text{where }c,b_j\in\N.
$$
Put $q_0=q_k=0$ and write the two sides of the $i$th relation as
$$
u_i=q_{i-1}+(d_i-d_{i-1})v
\quad\text{and}\quad t_i=q_i+a_iv.
$$
For the first relation, $\ell(u_1)=d_1\ge2$; for the last, $\ell(t_k)=a_k\ge2$. Every remaining side contains a $q_j$ and a positive multiple of $v$. Thus $\ell(u_i),\ell(t_i)\ge2$ for all $i$.

Two elements of $F$ have the same image in $V(A)$ if and only if one can be obtained from the other by finitely many replacements of $a+u_i$ by $a+t_i$ or vice versa, with $a\in F$ and $1\le i\le k$ at each step. Since $\ell(a+u_i),\ell(a+t_i)\ge2$ and $\ell(q_j)=1$, no replacement can start from $q_j$. Thus the class of $q_j$ is $\{q_j\}$, which contains neither $0$ nor $mv$ for any $m\in\N$. Consequently, $q_j\ne0$ and $q_j\ne mv$ in $V(A)$.
\end{proof}

The generators $q_i$ can be represented by explicit idempotents. For $1\le i<k$, let $U_i$ be the $d_i\times N_i$ block consisting of the first $N_i$ columns and first $d_i$ rows of $X$, and let $V_i$ be the $N_i\times d_i$ block consisting of the first $N_i$ rows and first $d_i$ columns of $Y$. The entries in the first $N_i$ columns of $X$ below row $d_i$ vanish, so $YX=I_n$ implies $V_iU_i=I_{N_i}$. Put
$$
\epsilon_i=I_{d_i}-U_iV_i\quad\text{and}\quad P_i=\epsilon_iA^{d_i}.
$$
Then $\epsilon_i^2=I_{d_i}-2U_iV_i+U_i(V_iU_i)V_i=\epsilon_i$, and $V_i\epsilon_i=0$, $\epsilon_iU_i=0$ and $U_iV_i+\epsilon_i=I_{d_i}$.
Let $\phi_i\colon A^{N_i}\oplus P_i\to A^{d_i}$ and $\psi_i\colon A^{d_i}\to A^{N_i}\oplus P_i$ be the maps defined by
$$
\phi_i(a,p)=U_ia+p\quad\text{and}\quad
\psi_i(z)=(V_iz,\epsilon_iz).
$$
Then
$$
\begin{aligned}
&\psi_i\phi_i(a,p)=(V_iU_ia+V_ip,\epsilon_iU_ia+\epsilon_ip)=(a,p),\\
&\phi_i\psi_i(z)=(U_iV_i+\epsilon_i)z=z.
\end{aligned}
$$
Thus
$$
A^{N_i}\oplus P_i\cong A^{d_i}\quad\text{and}\quad [P_i]=q_i.
$$
Preusser's isomorphism maps $q_i$ to $[\epsilon_i]=[P_i]$, as specified in \cite[proof of Theorem~14, Part~I, p.~136]{P}.

\section{The embedding criterion}
\label{sec:criterion}

For the algebra $A$ of an LV-rose, the relations in $V(A)$ from Proposition~\ref{prop:monoid} determine when a faithful integer-valued projective rank function exists. The same proposition shows that $A$ is right hereditary, so Lemma~\ref{lem:cohn} applies once such a rank function has been found.

Let $\rho$ be an integer-valued projective rank function on a unital ring $R$. For maps $\alpha\colon P\to Q$ between finitely generated projectives, let $r_\rho$ be the function defined by
$$
r_\rho(\alpha)=\min\{\rho(H):\alpha=\gamma\beta,
\ \beta\colon P\to H,\ \gamma\colon H\to Q,
\ H\text{ finitely generated projective}\}.
$$
The possible ranks $\rho(H)$ form a nonempty subset of $\N$, since one may take $H=P$. Cohn calls $r_\rho(\alpha)$ the \emph{inner $\rho$-rank} \cite[\S4.6, pp.~192--193]{C}. For square matrices put
$$
\Sigma_\rho=\mathop{\textstyle\bigcup}\limits_{m\ge1}\{M\in M_m(R):r_\rho(M)=m\}.
$$
A matrix in $\Sigma_\rho$ is called \emph{$\rho$-full}. For an idempotent $E\in M_s(R)$, put $P=ER^s$. The projection $E$ factors through $P$, so $r_\rho(E)\le\rho(P)$. If $E=\gamma\beta$ factors through $H$, then $E|_P=\id_P$ implies $\id_P=(E\gamma)(\beta|_P)$. Thus $P$ is a direct summand of $H$, whence $\rho(P)\le\rho(H)$. Since $\rho(P)\le\rho(H)$ for every finitely generated projective module $H$ through which $E$ factors, the definition of $r_\rho$ implies $\rho(P)\le r_\rho(E)$. Together with $r_\rho(E)\le\rho(P)$, this yields $r_\rho(E)=\rho(P)=\rho(ER^s)$.

\begin{lemma}
\label{lem:cohn}
Suppose that $R$ is right hereditary and $\rho$ is a faithful integer-valued projective rank function. Then the universal localization $\lambda\colon R\to R_{\Sigma_\rho}$ is injective, $R_{\Sigma_\rho}$ is a division ring, and $\rk_{R_{\Sigma_\rho}}\lambda(B)=r_\rho(B)$ for every rectangular matrix $B$ over $R$.
\end{lemma}

\begin{proof}
Let $\alpha\colon P\to Q$ and $\beta\colon Q\to T$ be maps between finitely generated projectives with $\beta\alpha=0$. Put $N=\ker\beta$ and $I=\im\beta$, and let $\beta_0\colon Q\to I$ be the map induced by $\beta$. Since $R$ is right hereditary, $I$ is projective. Choose $s\colon I\to Q$ with $\beta_0s=\id_I$. For $q\in Q$, $q=(q-s\beta_0(q))+s\beta_0(q)$ and $\beta_0(q-s\beta_0(q))=0$. Since $\ker\beta_0=N$, $Q=N+s(I)$. If $x\in I$ and $s(x)\in N$, then $x=\beta_0s(x)=0$. Hence $Q=N\oplus s(I)\cong N\oplus I$, so both $N$ and $I$ are finitely generated projective.

Let $\iota_N\colon N\to Q$ and $\iota_I\colon I\to T$ be the inclusions. Since $\beta\alpha=0$, let $\alpha_0\colon P\to N$ be the map defined by $\alpha_0(p)=\alpha(p)$. Then
$$
\alpha=\iota_N\alpha_0\quad\text{and}\quad \beta=\iota_I\beta_0.
$$
These factorizations imply
$$
r_\rho(\alpha)+r_\rho(\beta)
\le\rho(N)+\rho(I)=\rho(Q).
$$
This verifies Cohn's nullity condition for $\rho$ to be a Sylvester rank function \cite[\S4.6, p.~193]{C}. Since $\rho$ is integer-valued, Cohn's result on localization \cite[Proposition~4.6.13, pp.~194--195]{C} implies that $R_{\Sigma_\rho}$ is a division ring and $\rk_{R_{\Sigma_\rho}}\lambda(B)=r_\rho(B)$. Moreover,
$$
\ker\lambda=\bigl\langle e_{ij}:E=(e_{ij})\in M_s(R),\ E^2=E,
\ \rho(ER^s)=0\bigr\rangle.
$$
For an idempotent $E$ occurring here, $\rho(ER^s)=0$ and $\rho$ is faithful, so $ER^s=0$. If $\mathsf{e}_1,\ldots,\mathsf{e}_s$ are the standard columns of $R^s$, then
$$
E\mathsf{e}_j\in ER^s=0\quad\text{for }1\le j\le s,
\quad\text{hence}\quad
E=\begin{pmatrix}E\mathsf{e}_1&\cdots&E\mathsf{e}_s\end{pmatrix}=0.
$$
Therefore every generator of $\ker\lambda$ is zero, and $\ker\lambda=0$.
\end{proof}

\begin{theorem}
\label{thm:embedding}
Let $A=L_K(E,\omega)$ be an LV-algebra of a row-finite weighted graph. The following conditions are equivalent:
\begin{enumerate}
\renewcommand{\theenumi}{\roman{enumi}}
\renewcommand{\labelenumi}{\textup{(\theenumi)}}
\item $A$ embeds in a division ring.
\item $E$ has one vertex and, for the ordered weights of its $n$ structured loops, $w_n=n$ and $w_j>j$ for $1\le j<n$.
\item $E$ has one vertex and, with the distinct weights and multiplicities of Section~\ref{sec:preliminaries}, $d_k=N_k=n$ and $d_i>N_i$ for $1\le i<k$.
\end{enumerate}
When these conditions hold, the projective rank function on $A$ is unique. It is integer-valued and faithful, and is defined by $\rho(v)=1$ and $\rho(q_i)=d_i-N_i$ for $1\le i<k$.
\end{theorem}

\begin{proof}
Suppose first that $\lambda\colon A\hookrightarrow F$ is an embedding into a division ring. Each vertex of $E$ belongs to the normal form basis in Preusser \cite[Remark~7(e), p.~129]{P} and is thus nonzero in $A$. If $u,v$ are distinct vertices, then
$$
\lambda(u)\ne0,\qquad\lambda(v)\ne0\quad\text{and}\quad
\lambda(u)\lambda(v)=\lambda(uv)=0,
$$
which is impossible in $F$. Hence $E^0=\{v\}$. Since $E$ is row-finite, this rose is finite, with $1_A=v$. Moreover, $\lambda(v)^2=\lambda(v)\ne0$, hence $\lambda(v)=1_F$.

Let $\rho_F\colon V(A)\to\N$ be the function defined by $\rho_F([P])=\dim_F(P\otimes_AF)$ for each finitely generated projective $A$-module $P$. Since $(P\oplus Q)\otimes_AF\cong(P\otimes_AF)\oplus(Q\otimes_AF)$ and $A\otimes_AF\cong F$, we have
$$
\rho_F(P\oplus Q)=\dim_F\bigl((P\otimes_AF)\oplus(Q\otimes_AF)\bigr)=\rho_F(P)+\rho_F(Q)\quad\text{and}\quad\rho_F(A)=\dim_F F=1.
$$
Thus $\rho_F$ is an integer-valued projective rank function. To check that it is faithful, write $P\cong EA^s$ with $E^2=E\in M_s(A)$, and let $\mathsf{e}_j$ be the standard columns. Let $\Phi_E\colon(EA^s)\otimes_AF\to\lambda(E)F^s$ and $\Psi_E\colon\lambda(E)F^s\to(EA^s)\otimes_AF$ be the maps defined by
$$
\Phi_E(z\otimes f)=\lambda(z)f\quad\text{and}\quad
\Psi_E(\xi)=\sum_{j=1}^s E\mathsf{e}_j\otimes\xi_j.
$$
The maps $\Phi_E$ and $\Psi_E$ are mutually inverse.
Here $\lambda$ acts entrywise on columns. Indeed, $Ez=z$ and $\lambda(E)\xi=\xi$ imply
$$
\Phi_E\Psi_E(\xi)=\lambda(E)\xi=\xi\;\text{and}\;\Psi_E\Phi_E(z\otimes f)=\sum_jE\mathsf{e}_j\otimes\lambda(z_j)f=\sum_jE\mathsf{e}_jz_j\otimes f=Ez\otimes f=z\otimes f.
$$
If $P\ne0$, then $E\ne0$. Since $\lambda$ is injective, $\lambda(E)\ne0$. Therefore $P\otimes_AF\cong\lambda(E)F^s\ne0$ and $\rho_F(P)\ge1$.

Put $s_i=\rho_F(q_i)$ and $s_0=s_k=0$. Applying $\rho_F$ to the relations in Proposition~\ref{prop:monoid} and summing the resulting recurrence, we obtain
$$
\begin{aligned}
s_i-s_{i-1}&=d_i-d_{i-1}-a_i,\\
s_i&=\sum_{h=1}^i(d_h-d_{h-1}-a_h)
=d_i-\sum_{h=1}^i a_h=d_i-N_i.
\end{aligned}
$$
For $i<k$, $q_i\ne0$ and $\rho_F$ is faithful, so $d_i-N_i=s_i>0$. For $i=k$,
$$
0=s_k=d_k-N_k=d_k-n.
$$
Thus \textup{(iii)} holds.

Conversely, assume \textup{(iii)}. Let $\rho$ be the assignment on the generators defined by $\rho(v)=1$, $\rho(q_i)=d_i-N_i$ for $i<k$, and $\rho(q_0)=\rho(q_k)=0$. For each $1\le i\le k$,
$$
\begin{aligned}
\rho(q_{i-1})+d_i-d_{i-1}
&=d_{i-1}-N_{i-1}+d_i-d_{i-1}\\
&=d_i-N_i+a_i=\rho(q_i)+a_i.
\end{aligned}
$$
Since these are the relations presenting $V(A)$ in Proposition~\ref{prop:monoid}, the assignment induces an additive map $\rho\colon V(A)\to\N$ with $\rho([A])=1$. If $[P]=cv+\sum_{i=1}^{k-1}b_iq_i\ne0$, where $c,b_i\in\N$, then $c+\sum_i b_i\ge1$, and hence
$$
\rho(P)=c+\sum_{i=1}^{k-1}b_i(d_i-N_i)
\ge c+\sum_{i=1}^{k-1}b_i\ge1.
$$
Since $\rho$ is faithful and $A$ is right hereditary by Proposition~\ref{prop:monoid}, Lemma~\ref{lem:cohn} implies that $A\to A_{\Sigma_\rho}$ is injective and $A_{\Sigma_\rho}$ is a division ring.

For any projective rank function $\delta\colon V(A)\to\mathbb R_{\ge0}$, the same recurrence implies
$$
\delta(v)=1=\rho(v)
\quad\text{and}\quad
\delta(q_i)=\sum_{h=1}^i(d_h-d_{h-1}-a_h)=d_i-N_i=\rho(q_i).
$$
Consequently, for every $x=cv+\sum_i b_iq_i\in V(A)$, we have
$$
\delta(x)=c+\sum_i b_i\delta(q_i)
=c+\sum_i b_i\rho(q_i)=\rho(x).
$$
Thus $\delta=\rho$.

It remains to compare \textup{(ii)} and \textup{(iii)}. Assume \textup{(iii)} and let $w_j=d_i$, with $j<n$. Then
$$
w_j-j=\begin{cases}
d_i-j\ge d_i-N_i>0,&i<k,\\
n-j>0,&i=k.
\end{cases}
$$
Also $w_n=d_k=n$, so \textup{(ii)} holds. Conversely, under \textup{(ii)}, for $i<k$ we have $N_i<n$ and
$$
d_i-N_i=w_{N_i}-N_i>0,
\quad\text{while}\quad d_k=w_n=n.
$$
Hence \textup{(iii)} holds.
\end{proof}

\begin{corollary}
\label{cor:stablefree}
Suppose that the LV-algebra $A$ embeds in a division ring, and let $\rho$ be its projective rank function from Theorem~\ref{thm:embedding}. Then every finitely generated projective $A$-module is stably free. More precisely, if $[P]=cv+\sum_{i=1}^{k-1}b_iq_i$, with $c,b_i\in\N$, then $P\oplus A^s\cong A^{\rho(P)+s}$, where $s=\sum_{i=1}^{k-1}b_iN_i$.
\end{corollary}

\begin{proof}
The expression for $[P]$ means $P\cong A^c\oplus\bigoplus_iP_i^{b_i}$. Using the isomorphisms $P_i\oplus A^{N_i}\cong A^{d_i}$ constructed in Section~\ref{sec:preliminaries}, we obtain
$$
\begin{aligned}
P\oplus A^s
&\cong A^c\oplus\bigoplus_i(P_i\oplus A^{N_i})^{b_i}\\
&\cong A^{c+\sum_i b_id_i}.
\end{aligned}
$$
The exponent is
$$
c+\sum_i b_id_i
=c+\sum_i b_i(d_i-N_i)+\sum_i b_iN_i
=\rho(P)+s,
$$
so $P\oplus A^s\cong A^{\rho(P)+s}$.
\end{proof}

Following Cohn \cite[Theorem~1.6.1, p.~35]{C}, we call a nonzero unital ring a \emph{semifir} if every finitely generated left or right ideal is free of unique rank.

For a square matrix $M\in M_m(R)$ over a unital ring $R$, the \emph{inner rank} of $M$ is the least $q\ge0$ such that $M=BC$ with $B\in M_{m\times q}(R)$ and $C\in M_{q\times m}(R)$. The matrix $M$ is \emph{full} if its inner rank is $m$ \cite[\S4.5, p.~179]{C}.

\begin{corollary}
\label{cor:semifir}
An LV-algebra $A=L_K(E,\omega)$ is a semifir if and only if $E$ is a rose with $n$ structured loops, all of weight $n$. In this case $\Sigma_\rho$ is the set of all full square matrices, and the resulting localization is the one in Cohn's embedding theorem for semifirs \cite[Corollary~4.5.9, p.~182]{C}.
\end{corollary}

\begin{proof}
Cohn proves that a semifir embeds in a division ring \cite[Corollary~4.5.9, p.~182]{C}. The embedding criterion in Theorem~\ref{thm:embedding} then implies that $E$ is a rose. We first show that every finitely generated projective right module is free. Such a module is a direct summand of some $A^r$, so it suffices to induct on $r$ for a finitely generated projective $P\subseteq A^r$. For $r=0$, $P=0$; for $r=1$, $P$ is a finitely generated right ideal and hence free.

For $r>1$, let $\pi\colon P\to A$ be the map defined by $\pi(x)=x_r$ for $x=(x_1,\ldots,x_r)^T\in P$, and put $I=\pi(P)$. Since $I$ is a finitely generated right ideal, choose a free basis $a_1,\ldots,a_t$ of $I$ and elements $p_1,\ldots,p_t\in P$ with $\pi(p_j)=a_j$. Let $\sigma\colon I\to P$ be the map defined by $\sigma(\sum_{j=1}^t a_jc_j)=\sum_{j=1}^t p_jc_j$ for $c_j\in A$. Then $\pi\sigma(x)=x$ for $x\in I$. Put $N=\ker\pi$. For every $p\in P$, $p=(p-\sigma\pi(p))+\sigma\pi(p)$ and $\pi(p-\sigma\pi(p))=0$. If $x\in I$ and $\sigma(x)\in N$, then $x=\pi\sigma(x)=0$. Hence $P=N\oplus\sigma(I)\cong N\oplus I$. The module $N$ is finitely generated and projective, and $N\subseteq A^{r-1}$. By induction, $N\cong A^s$ for some $s$, and consequently
$$
P\cong N\oplus I\cong A^s\oplus A^t\cong A^{s+t}.
$$

In particular, $P_i\cong A^{m_i}$ for some $m_i\in\N$, so $q_i=[P_i]=m_iv$. Proposition~\ref{prop:monoid} excludes $q_i=m_iv$ for $i<k$, so $k=1$. The equality $d_k=n$ in Theorem~\ref{thm:embedding} thus reduces to $d_1=n$.

Conversely, if all $n$ weights equal $n$, the presentation of $V(A)$ in Proposition~\ref{prop:monoid} reduces to
$$
V(A)\cong\langle v\mid nv=nv\rangle\cong\N,
$$
where the second isomorphism sends $v$ to $1$. Thus every finitely generated projective right module $P$ satisfies $[P]=mv$ for a unique $m$, and hence $P\cong A^m$. In particular,
$$
A^m\cong A^{m'}\ \Longrightarrow\ mv=m'v\ \Longrightarrow\ m=m'.
$$
Since $A$ is right hereditary, every finitely generated right ideal is projective and therefore free of unique rank.

For a finitely generated projective left module $L$, put $L^\vee=\operatorname{Hom}_{A^{\mathrm{op}}}(L,A)$. Since $L$ is finitely generated and projective, the map $L\to\operatorname{Hom}_A(L^\vee,A)$ defined by $\ell\mapsto(f\mapsto f(\ell))$ is an isomorphism. Thus
$$
L\cong\operatorname{Hom}_A(L^\vee,A)
\cong\operatorname{Hom}_A(A^m,A)
\cong {}_AA^m
$$
for the unique $m$ with $L^\vee\cong A^m$. Since $A$ is also left hereditary, every finitely generated left ideal is free of unique rank. Thus $A$ is a semifir.

Finally, every finitely generated projective module $H$ is isomorphic to $A^t$ with $\rho(H)=t$. For $M\in M_m(A)$, the definition of $r_\rho$ therefore becomes
$$
r_\rho(M)=\min\{t:M=BC,\ B\in M_{m\times t}(A),\ C\in M_{t\times m}(A)\}.
$$
Thus $r_\rho(M)$ equals the inner rank of $M$, and $\Sigma_\rho$ consists exactly of the full square matrices.
\end{proof}

\begin{example}
\label{ex:domain-no-embedding}
Consider the LV-rose with weights $(2,2,4,4)$. As in Hazrat's diagram \cite[Example~5.5, p.~876]{Hazrat2013GradedStructure}, we label each structured loop by its corresponding edges.
\begin{center}
\begin{tikzpicture}[>={Stealth[length=1.3mm,width=.9mm]}]
    \node[inner sep=1pt] at (0,0) (v) {$\bullet$};
    \node at (0,-0.7) {$v$};
    \draw[->] (v) .. controls (-1.9,0.3) and (-0.3,1.9) .. node[above left] {$e_1,e_2$} (v);
    \draw[->] (v) .. controls (0.3,1.9) and (1.9,0.3) .. node[above right] {$f_1,f_2$} (v);
    \draw[->] (v) .. controls (-0.3,-1.9) and (-1.9,-0.3) .. node[below left] {$g_1,\ldots,g_4$} (v);
    \draw[->] (v) .. controls (1.9,-0.3) and (0.3,-1.9) .. node[below right] {$h_1,\ldots,h_4$} (v);
\end{tikzpicture}
\end{center}
For these weights, the presentation in Proposition~\ref{prop:monoid} reads
$$
V(A)\cong\langle v,q\mid2v=q+2v,\ q+2v=2v\rangle.
$$
Here $q\ne0$ in $V(A)$. The relation $2v=q+2v$ implies $\iota(q)=0$ in $K_0(A)$ and $\rho(q)=0$ for every projective rank function $\rho$. Suppose that $\lambda\colon A\hookrightarrow F$ is an embedding into a division ring. Let $\rho_F\colon V(A)\to\mathbb R_{\ge0}$ be the projective rank function defined by $\rho_F([Q])=\dim_F(Q\otimes_AF)$. Choose an idempotent $E\in M_s(A)$ such that $[EA^s]=q$. Since $q\ne0$, we have $E\ne0$, and injectivity of $\lambda$ implies $\lambda(E)\ne0$. Hence
$$
0=\rho_F(q)=\dim_F((EA^s)\otimes_AF)=\dim_F\bigl(\lambda(E)F^s\bigr)>0,
$$
a contradiction. Thus $A$ does not embed in a division ring. Since the graph has one vertex, the criterion of Hazrat and Preusser \cite[Theorem~41, p.~1079]{HP} implies that $A$ is a domain. In particular, the domain condition alone is insufficient.
\end{example}

\section{The universal division ring of fractions}
\label{sec:construction}

For the rest of the paper, fix an LV-algebra $A=L_K(E,\omega)$ satisfying the embedding criterion of Theorem~\ref{thm:embedding}, with its unique projective rank function $\rho$, unless otherwise specified. Since $A$ is right hereditary by Proposition~\ref{prop:monoid} and $\rho$ is faithful and integer-valued, Lemma~\ref{lem:cohn} identifies $A_{\Sigma_\rho}$ as a division ring containing $A$.

We describe this localization by adjoining inverses to the matrices in $\Sigma_\rho$, as in Cohn \cite[Theorem~4.1.3, p.~156]{C}. For each $M=(m_{ij})\in\Sigma_\rho\cap M_m(A)$, we introduce an $m\times m$ matrix of new symbols $Z_M=(z^M_{ij})$. Let $\mathcal J$ be the two-sided ideal generated by the entries of $MZ_M-I_m$ and $Z_MM-I_m$ for all such $M$, and set
$$
D_K(E,\omega)=\bigl(A*_K K\langle z^M_{ij}:M\in\Sigma_\rho\cap M_m(A),\ m\ge1,\ 1\le i,j\le m\rangle\bigr)\big/\mathcal J.
$$
The free product is taken over $K$, so $cz^M_{ij}=z^M_{ij}c$ for $c\in K$. The equations $MZ_M=Z_MM=I_m$ can be written as
$$
\sum_{j=1}^m m_{ij}z^M_{j\ell}=\delta_{i\ell}
\quad\text{and}\quad
\sum_{j=1}^m z^M_{ij}m_{j\ell}=\delta_{i\ell}
\quad\text{for }1\le i,\ell\le m.
$$
Let $f\colon A\to B$ be a unital homomorphism such that every $f(M)$, $M\in\Sigma_\rho$, is invertible. For $c\in K$ and $M\in\Sigma_\rho\cap M_m(A)$, the entries of $M$ commute with $c$, so
$$
(f(c)I_m)f(M)=f(M)(f(c)I_m),
\quad\text{hence}\quad
(f(c)I_m)f(M)^{-1}=f(M)^{-1}(f(c)I_m).
$$
Thus the entries of $f(M)^{-1}$ commute with $f(K)$. Let $\widetilde f\colon D_K(E,\omega)\to B$ be the homomorphism defined by
$$
\widetilde f(a)=f(a)\quad\text{for }a\in A\quad\text{and}\quad
\widetilde f(Z_M)=f(M)^{-1}\quad\text{for }M\in\Sigma_\rho.
$$
Indeed, these substitutions respect the relations of $D_K(E,\omega)$, since $f(M)f(M)^{-1}=I_m=f(M)^{-1}f(M)$ for $M\in\Sigma_\rho\cap M_m(A)$. Conversely every extension of $f$ must have these values, so it is unique. Thus $D_K(E,\omega)=A_{\Sigma_\rho}$. We henceforth identify $A$ with its image under the canonical embedding $\lambda\colon A\hookrightarrow D_K(E,\omega)$.

An \emph{epic $A$-field} is a division ring $F$ with a unital homomorphism $f\colon A\to F$ whose image generates $F$ as a division ring. The homomorphism $f$ need not be injective. Its \emph{singular kernel} is
$$
\mathcal P_f=\mathop{\textstyle\bigcup}\limits_{m\ge1}\{M\in M_m(A):f(M)\text{ is singular over }F\}.
$$
For two epic $A$-fields $(F_1,f_1)$ and $(F_2,f_2)$, a \emph{specialization} from $(F_1,f_1)$ to $(F_2,f_2)$ is represented by a \emph{local} homomorphism $h\colon S\to F_2$ from a subring $f_1(A)\subseteq S\subseteq F_1$, with $h\circ f_1=f_2$. Here $h$ is local if $h(s)\ne0$ implies $s^{-1}\in S$ for every $s\in S$. Two such homomorphisms represent the same specialization if and only if they agree on a common subring containing $f_1(A)$ and the common restriction is local \cite[\S4.1, p.~154]{C}.

According to Cohn's specialization theorem \cite[Theorem~4.4.1, p.~172]{C}, there is such a specialization if and only if $\mathcal P_{f_1}\subseteq\mathcal P_{f_2}$. Following Cohn \cite[\S4.4, p.~176]{C}, we call an epic $A$-field $(F,f)$ a \emph{universal $A$-field} if it specializes to every epic $A$-field. If $f$ is injective, then $F$ is a \emph{universal division ring of fractions} of $A$.

\begin{theorem}
\label{thm:universal}
The division ring $D_K(E,\omega)$ constructed above is a universal division ring of fractions of $A$. It is unique up to a unique $A$-isomorphism. For a unital homomorphism $f\colon A\to F$ into a division ring, a homomorphism $D_K(E,\omega)\to F$ extending $f$ exists if and only if $f$ inverts every matrix in $\Sigma_\rho$. When it exists, it is unique and injective.
\end{theorem}

\begin{proof}
Let $F_0$ be the division subring generated by $A$ in $D_K(E,\omega)$. For $M\in\Sigma_\rho\cap M_m(A)$, choose $U,V\in\operatorname{GL}_m(F_0)$ and $0\le r\le m$ by row and column reduction over $F_0$, so that $UMV=\operatorname{diag}(I_r,0)$. Since $U,M,V$ are invertible over $D_K(E,\omega)$, we get
$$
r=\rk_{D_K(E,\omega)}\operatorname{diag}(I_r,0)=\rk_{D_K(E,\omega)}(UMV)=m.
$$
Thus $UMV=I_m$ and $M^{-1}=VU\in M_m(F_0)$. Since $F_0$ contains $A$ and all entries of the adjoined inverses, the presentation implies $F_0=D_K(E,\omega)$. In particular, $D_K(E,\omega)$ is an epic $A$-field.

Let $f\colon A\to F$ define any epic $A$-field. Put $P_F=P\otimes_AF$ for each finitely generated projective $A$-module $P$, and let $\delta_f\colon V(A)\to\mathbb R_{\ge0}$ be the function defined by $\delta_f([P])=\dim_F P_F$. The isomorphisms
$$
(P\oplus Q)_F\cong P_F\oplus Q_F
\quad\text{and}\quad A_F\cong F
$$
imply
$$
\delta_f(P\oplus Q)=\delta_f(P)+\delta_f(Q)
\quad\text{and}\quad\delta_f(A)=1.
$$
Thus $\delta_f$ is a projective rank function. The uniqueness assertion of Theorem~\ref{thm:embedding} implies $\delta_f=\rho$. Suppose that $f(M)$ is invertible for some $M\in M_m(A)$. A factorization $M=\gamma\beta$ through a finitely generated projective $H$ induces
$$
F^m\xrightarrow{\ \beta\otimes1\ }H_F
\xrightarrow{\ \gamma\otimes1\ }F^m,
\quad\text{with}\quad
f(M)=(\gamma\otimes1)(\beta\otimes1).
$$
Consequently
$$
m=\rk_F f(M)\le\dim_F H_F=\delta_f(H)=\rho(H).
$$
For every such $H$, $m\le\rho(H)$, so $m\le r_\rho(M)$. Also, $r_\rho(M)\le m$ because $M$ factors through $A^m$.
Hence $M\in\Sigma_\rho$. In terms of singular kernels, this proves
$$
M\notin\mathcal P_f\ \Longrightarrow\ M\notin\mathcal P_\lambda,
\quad\text{and therefore}\quad
\mathcal P_\lambda\subseteq\mathcal P_f.
$$
It follows from Cohn's specialization theorem \cite[Theorem~4.4.1, p.~172]{C} that there is a specialization $D_K(E,\omega)\to F$. It is represented by a local homomorphism $h\colon S\to F$ with $A\subseteq S\subseteq D_K(E,\omega)$ and $h|_A=f$. The unlabeled arrows below are inclusions.
\begin{center}
\begin{tikzpicture}[>={Stealth[length=1.3mm,width=.9mm]}]
    \node (A) at (0,0) {$A$};
    \node (S) at (2,0) {$S$};
    \node (D) at (4,0) {$D_K(E,\omega)$};
    \node (F) at (2,-1.2) {$F$};
    \draw[->] (A) -- (S);
    \draw[->] (S) -- (D);
    \draw[->] (S) -- node[right] {$h$} (F);
    \draw[->] (A) -- node[below left] {$f$} (F);
\end{tikzpicture}
\end{center}
Hence $D_K(E,\omega)$ is universal as an $A$-field. In Cohn's terminology, $\mathcal P_\lambda$ is the least prime matrix ideal of $A$.

If $(D',\lambda')$ is another universal division ring of fractions of $A$, the specializations $D_K(E,\omega)\to D'$ and $D'\to D_K(E,\omega)$ imply $\mathcal P_\lambda\subseteq\mathcal P_{\lambda'}$ and $\mathcal P_{\lambda'}\subseteq\mathcal P_\lambda$. By the last assertion of Cohn's theorem \cite[Theorem~4.4.1, p.~172]{C}, equality of these singular kernels implies that there is an $A$-isomorphism $D_K(E,\omega)\to D'$. To check that it is unique, let $\theta_1,\theta_2\colon D_K(E,\omega)\to D'$ be $A$-homomorphisms. They agree on $A$, and the set $\{a\in D_K(E,\omega):\theta_1(a)=\theta_2(a)\}$ is closed under addition, multiplication and inverses because
$$
\theta_i(a+b)=\theta_i(a)+\theta_i(b),\quad \theta_i(ab)=\theta_i(a)\theta_i(b)\quad\text{and}\quad \theta_i(a^{-1})=\theta_i(a)^{-1}\quad\text{whenever }a\ne0.
$$
Since $D_K(E,\omega)$ is the division ring generated by $A$, it follows that $\theta_1=\theta_2$.

Finally, if $g\colon D_K(E,\omega)\to F$ extends a unital homomorphism $f\colon A\to F$, then
$$
f(M)g(M^{-1})=I_m=g(M^{-1})f(M)\quad\text{for }M\in\Sigma_\rho,
$$
so $f$ inverts $\Sigma_\rho$. Conversely, whenever $f$ inverts $\Sigma_\rho$, the relations of $D_K(E,\omega)$ hold after replacing $a\in A$ by $f(a)$ and $Z_M$ by $f(M)^{-1}$. Thus $f$ has a unique extension $g\colon D_K(E,\omega)\to F$ defined by
$$
g(a)=f(a)\quad\text{for }a\in A
\quad\text{and}\quad
\bigl(g(z^M_{ij})\bigr)_{ij}=f(M)^{-1}\quad\text{for }M\in\Sigma_\rho.
$$ Since $g(1)=1$, its kernel is a proper ideal of the division ring $D_K(E,\omega)$; therefore $\ker g=0$.
\end{proof}

\begin{remark}
Let $\theta\colon A\to A'$ be induced by a weighted graph isomorphism. The projective rank function on $A$ is unique, so $\rho'(P\otimes_AA')=\rho(P)$. Applying $\theta$ and $\theta^{-1}$ to factorizations through finitely generated projective modules shows that $r_{\rho'}(\theta(M))=r_\rho(M)$. Therefore $\theta(\Sigma_\rho)=\Sigma_{\rho'}$, so there is a unique $K$-isomorphism $D_K(E,\omega)\to D_K(E',\omega')$ extending $\theta$ and defined by
$$
a\longmapsto\theta(a)\quad\text{for }a\in A
\quad\text{and}\quad
M^{-1}\longmapsto\theta(M)^{-1}\quad\text{for }M\in\Sigma_\rho.
$$ For a rose, this applies in particular to any relabeling of loops with the same weights.
\end{remark}

\section{Free generators and classification}
\label{sec:free-coordinates}

The equality $w_n=n$ in Theorem~\ref{thm:embedding} makes the matrices $X$ and $Y$ of Section~\ref{sec:preliminaries} square, with $XY=YX=I_n$. We factor $X$ to find generators of $D_K(E,\omega)$ over $K$ and use the specialization property in Theorem~\ref{thm:universal} to prove that they are free. Put $\mathcal I=\{(i,j):1\le j<i\le w_j\}$.

\begin{lemma}
\label{lem:gauss}
\label{lem:gauss-factors}
Every leading principal submatrix $X[q]=(x_{ij})_{1\le i,j\le q}$ is invertible over $D_K(E,\omega)$. There is a unique factorization $X=LPU$, where $L=(l_{ij})$ is lower unitriangular, $U=(u_{ij})$ is upper unitriangular and $P=\operatorname{diag}(p_1,\ldots,p_n)$ has nonzero diagonal entries. If $T=(t_{ij})=U^{-1}$, then $l_{ij}=0$ and $t_{ji}=0$ whenever $i>w_j$. The elements $p_j$ for $1\le j\le n$ and $l_{ij},t_{ji}$ for $(i,j)\in\mathcal I$, together with $K$, generate $D_K(E,\omega)$ as a division ring.
\end{lemma}

\begin{proof}
Since $w_j\ge j$ for $1\le j\le n$ by Theorem~\ref{thm:embedding}, define a unital $K$-algebra homomorphism $\varepsilon\colon A\to K$ by
$$
\varepsilon(x_{ij})=\delta_{ij}
\quad\text{and}\quad
\varepsilon(y_{ji})=\delta_{ij}.
$$
Indeed, the assignments respect the prescribed zeros and satisfy $\varepsilon(X)\varepsilon(Y)=\varepsilon(Y)\varepsilon(X)=I_n$. The function $V(A)\to\mathbb R_{\ge0}$ defined by $[H]\mapsto\dim_K(H\otimes_{A,\varepsilon}K)$ is a projective rank function. The uniqueness in Theorem~\ref{thm:embedding} identifies this function with $\rho$. If $X[q]$ factors through a finitely generated projective $H$, then $\varepsilon(X[q])$ factors through $H\otimes_{A,\varepsilon}K$, so
$$
q=\rk_K\varepsilon(X[q])
\le\dim_K(H\otimes_{A,\varepsilon}K)
=\rho(H).
$$
The factorization through $A^q$ implies $r_\rho(X[q])\le q$, and hence $r_\rho(X[q])=q$. Applying Cohn's rank formula \cite[Proposition~4.6.13, pp.~194--195]{C} as in Lemma~\ref{lem:cohn}, we get $\rk_{D_K(E,\omega)}(X[q])=q$. Hence $X[q]$ is invertible over $D_K(E,\omega)$.

We construct the factors by induction. The first pivot is $p_1=x_{11}\ne0$. For $q\ge2$, put $B=X[q-1]$ and write
$$
X[q]=\begin{pmatrix}B&c\\a&x_{qq}\end{pmatrix}.
$$
Set $p_q=x_{qq}-aB^{-1}c$. Direct multiplication shows that
$$
\begin{pmatrix}I&0\\-aB^{-1}&1\end{pmatrix}
X[q]
\begin{pmatrix}I&-B^{-1}c\\0&1\end{pmatrix}
=\begin{pmatrix}B&0\\0&p_q\end{pmatrix}.
$$
Therefore $p_q\ne0$. If $X[q-1]=L_{q-1}P_{q-1}U_{q-1}$, set
$$
\ell=aU_{q-1}^{-1}P_{q-1}^{-1}
\quad\text{and}\quad
u=P_{q-1}^{-1}L_{q-1}^{-1}c.
$$
The identities
$$
L_{q-1}P_{q-1}u=c,\qquad
\ell P_{q-1}U_{q-1}=a
\quad\text{and}\quad
\ell P_{q-1}u+p_q=aB^{-1}c+p_q=x_{qq}
$$
verify the block multiplication
$$
\begin{pmatrix}L_{q-1}&0\\\ell&1\end{pmatrix}
\begin{pmatrix}P_{q-1}&0\\0&p_q\end{pmatrix}
\begin{pmatrix}U_{q-1}&u\\0&1\end{pmatrix}
=X[q].
$$
This constructs $L,P,U$. Their entries are recovered recursively by $p_j=x_{jj}-\sum_{k<j}l_{jk}p_ku_{kj}$ for $1\le j\le n$, and by $l_{ij}=\bigl(x_{ij}-\sum_{k<j}l_{ik}p_ku_{kj}\bigr)p_j^{-1}$ and $u_{ji}=p_j^{-1}\bigl(x_{ji}-\sum_{k<j}l_{jk}p_ku_{ki}\bigr)$ for $j<i\le n$. If $LPU=L'P'U'$ is another such factorization, then
$$
L'^{-1}L=P'U'U^{-1}P^{-1}.
$$
The matrix $H=L'^{-1}L=P'U'U^{-1}P^{-1}$ is lower unitriangular and upper triangular, so
$$
H_{ij}=0\quad\text{for }i\ne j
\quad\text{and}\quad
H_{jj}=1=p'_jp_j^{-1}\quad\text{for }1\le j\le n.
$$
Hence $H=I_n$, $L'=L$, $P'=P$ and $U'=U$.

For $i>w_j$, induction on $j$ in the formula for $l_{ij}$ proves $l_{ij}=0$. At stage $j$, the inequalities $k<j$ and $i>w_j\ge w_k$ imply $l_{ik}=0$ by the preceding stages. Since $x_{ij}=0$,
$$
l_{ij}=\bigl(0-\sum_{k<j}0\cdot p_ku_{kj}\bigr)p_j^{-1}=0.
$$ Since $Y=X^{-1}$, we have
$$
Y=TP^{-1}L^{-1}
\quad\text{and}\quad
T=YLP.
$$
For $i>w_j$,
$$
t_{ji}=\sum_{k=1}^n y_{jk}l_{ki}p_i
=\sum_{k\le w_j}y_{jk}l_{ki}p_i=0.
$$
Here $y_{jk}=0$ for $k>w_j$, while $k\le w_j<i$ implies $l_{ki}=0$ because $L$ is lower triangular.

Let $D_0$ be the division subring of $D_K(E,\omega)$ generated by $K$, $p_j$ for $1\le j\le n$, and $l_{ij},t_{ji}$ for $(i,j)\in\mathcal I$. The finite inverse formula for unitriangular matrices implies
$$
X=LPT^{-1}\in M_n(D_0)
\quad\text{and}\quad
Y=TP^{-1}L^{-1}\in M_n(D_0).
$$
Thus $A\subseteq D_0$. Theorem~\ref{thm:universal} identifies $D_K(E,\omega)$ as a universal division ring of fractions of $A$. Since $A$ generates it as a division ring, $D_0=D_K(E,\omega)$.
\end{proof}

Put $r(E,\omega)=n+2|\mathcal I|=2\sum_{j=1}^n w_j-n^2$.

\begin{theorem}
\label{thm:free-coordinates}
The elements $p_j$ for $1\le j\le n$ and $l_{ij},t_{ji}$ for $(i,j)\in\mathcal I$ in Lemma~\ref{lem:gauss-factors} form a free generating set of $D_K(E,\omega)$ over $K$. In particular, $D_K(E,\omega)\cong_K\mathcal F_K(r(E,\omega))$.
\end{theorem}

\begin{proof}
Let $B$ be the free $K$-algebra on the symbols $\pi_j$ for $1\le j\le n$ and $\lambda_{ij},\theta_{ji}$ for $(i,j)\in\mathcal I$. According to Cohn \cite[\S5.4, p.~224]{C}, $B$ is a fir. It follows from his embedding theorem for semifirs \cite[Corollary~4.5.9, p.~182]{C} that $B$ has a universal division ring of fractions $F=U(B)$. In $F$, form a lower unitriangular matrix $L_0$ and an upper unitriangular matrix $T_0$ by setting
$$
(L_0)_{ij}=\lambda_{ij}
\quad\text{and}\quad
(T_0)_{ji}=\theta_{ji}
\qquad\text{for }(i,j)\in\mathcal I,
$$
with all other strictly triangular entries zero. Put
$$
P_0=\operatorname{diag}(\pi_1,\ldots,\pi_n),\qquad
X_0=L_0P_0T_0^{-1}
\quad\text{and}\quad
Y_0=T_0P_0^{-1}L_0^{-1}.
$$
The entries $\pi_j$ are nonzero. The unitriangular factors are invertible by the finite identity
$$
V^{-1}=\sum_{a=0}^{n-1}(-1)^a(V-I_n)^a.
$$
In particular,
$$
X_0Y_0=L_0P_0(T_0^{-1}T_0)P_0^{-1}L_0^{-1}=I_n\quad\text{and}\quad Y_0X_0=T_0P_0^{-1}(L_0^{-1}L_0)P_0T_0^{-1}=I_n.
$$

Since $T_0^{-1}$ is upper triangular, for $i>w_j$ we have
$$
(X_0)_{ij}
=\sum_{k\le j}(L_0)_{ik}\pi_k(T_0^{-1})_{kj}=0,
$$
since $i>w_j\ge w_k$ for every $k\le j$. Also
$$
(Y_0)_{ji}
=\sum_{k\ge i}(T_0)_{jk}\pi_k^{-1}(L_0^{-1})_{ki}=0,
$$
because $k\ge i>w_j$. Let $j\colon A\to F$ be the unital $K$-algebra homomorphism defined by
$$
j(X)=X_0\quad\text{and}\quad j(Y)=Y_0.
$$

For every $q$, the leading principal submatrix has the factorization
$$
X_0[q]=L_0[q]P_0[q](T_0[q])^{-1}.
$$
Its three factors are invertible. Let $F_0$ be the division subring of $F$ generated by $j(A)$. Applying the recursions in the proof of Lemma~\ref{lem:gauss-factors} to $X_0\in M_n(F_0)$ shows that $L_0,P_0,U_0\in M_n(F_0)$, where $U_0=T_0^{-1}$. Thus $T_0=U_0^{-1}\in M_n(F_0)$ by the finite inverse formula. Consequently
$$
\pi_j\in F_0\quad\text{for }1\le j\le n,\qquad
\lambda_{ij},\theta_{ji}\in F_0\quad\text{for }(i,j)\in\mathcal I.
$$
Thus $B\subseteq F_0\subseteq F$, and $F_0=F$ because $F=U(B)$ is generated by $B$ as a division ring. Hence $(F,j)$ is an epic $A$-field.

It follows from the universality assertion of Theorem~\ref{thm:universal} that $D_K(E,\omega)$ specializes to this epic $A$-field. Represent this specialization by a local homomorphism
$$
h:S\longrightarrow F,
\qquad A\subseteq S\subseteq D_K(E,\omega),
\qquad h|_A=j.
$$
We verify that $p_j\in S$ for $1\le j\le n$ and $l_{ij},t_{ji}\in S$ for $(i,j)\in\mathcal I$. At the first stage, $p_1=x_{11}\in S$ and $h(p_1)=\pi_1\ne0$. Since $h$ is local, $p_1^{-1}\in S$. The entries in the first column of $L$ and the first row of $U$ therefore belong to $S$ and map to the corresponding entries of $L_0,U_0$.

Suppose the pivots and the entries of $L,U$ constructed before stage $j$ belong to $S$ and map under $h$ to the corresponding entries of $P_0,L_0,U_0$. The pivot formula implies
$$
p_j=x_{jj}-\sum_{k<j}l_{jk}p_ku_{kj}\in S
\quad\text{and}\quad
h(p_j)=\pi_j\ne0.
$$
Hence $p_j^{-1}\in S$. The recursions imply $l_{ij},u_{ji}\in S$. Applying $h$ to the recursions proves by induction that
$$
h(p_j)=\pi_j,\qquad
h(l_{ij})=(L_0)_{ij}
\quad\text{and}\quad
h(u_{ji})=(U_0)_{ji}.
$$
Since $T=U^{-1}$ is a finite polynomial in $U-I_n$, its entries also belong to $S$ and $h(T)=T_0$.

Let $\beta\colon B\to D_K(E,\omega)$ be the $K$-algebra homomorphism defined by
$$
\beta(\pi_j)=p_j,\qquad
\beta(\lambda_{ij})=l_{ij}
\quad\text{and}\quad
\beta(\theta_{ji})=t_{ji}.
$$
Write $\iota_B:B\hookrightarrow F$ for the canonical embedding. The coordinate identities imply
$$
\beta(B)\subseteq S,\qquad h\circ\beta=\iota_B
\quad\text{and}\quad
\ker\beta\subseteq\ker(h\circ\beta)=\ker\iota_B=0.
$$
Hence $\beta$ is injective. The image $\beta(B)$ contains $K$ and the generators from Lemma~\ref{lem:gauss-factors}, so it generates $D_K(E,\omega)$ as a division ring. Hence $(D_K(E,\omega),\beta)$ is an epic $B$-field. The map $h$ represents a $B$-specialization to $F$, and the universality of $F=U(B)$ supplies a specialization in the reverse direction. The two specializations imply
$$
\mathcal P_\beta\subseteq\mathcal P_{\iota_B}
\quad\text{and}\quad
\mathcal P_{\iota_B}\subseteq\mathcal P_\beta.
$$
Since the singular kernels are equal, Cohn's theorem \cite[Theorem~4.4.1, p.~172]{C} implies that there is a $B$-isomorphism $\Phi:F\to D_K(E,\omega)$ with $\Phi\circ\iota_B=\beta$. It is unique because $F$ is generated by $\iota_B(B)$ as a division ring. Moreover,
$$
\Phi(j(X))=LPT^{-1}=X
\quad\text{and}\quad
\Phi(j(Y))=TP^{-1}L^{-1}=Y.
$$
Finally,
$$
n+2|\mathcal I|
=n+2\sum_{j=1}^n(w_j-j)
=2\sum_{j=1}^n w_j-n^2.
$$
Hence $D_K(E,\omega)\cong_K\mathcal F_K(r(E,\omega))$.
\end{proof}

\begin{corollary}
\label{cor:field-classification}
For two LV-roses whose algebras over the same field $K$ satisfy the embedding criterion of Theorem~\ref{thm:embedding},
$$
D_K(E,\omega)\cong_KD_K(E',\omega')
\quad\Longleftrightarrow\quad
r(E,\omega)=r(E',\omega').
$$
Furthermore, $Z(D_K(E,\omega))=K$.
\end{corollary}

\begin{proof}
Put $r=r(E,\omega)$ and $r'=r(E',\omega')$. The free generating sets in Theorem~\ref{thm:free-coordinates} determine $K$-isomorphisms
$$
\Phi:\mathcal F_K(r)\longrightarrow D_K(E,\omega)
\quad\text{and}\quad
\Phi':\mathcal F_K(r')\longrightarrow D_K(E',\omega').
$$
If $r=r'$, the composite $\Phi'\Phi^{-1}:D_K(E,\omega)\to D_K(E',\omega')$ is a $K$-isomorphism.

Conversely, suppose that $\psi:D_K(E,\omega)\to D_K(E',\omega')$ is a $K$-isomorphism. Let $\Theta\colon\mathcal F_K(r)\to\mathcal F_K(r')$ be the map defined by $\Theta=(\Phi')^{-1}\psi\Phi$. Then $\Theta$ is a $K$-isomorphism. If $s_1,\ldots,s_r$ and $s'_1,\ldots,s'_{r'}$ are the standard free generating sets of these free fields, then $\{\Theta(s_1),\ldots,\Theta(s_r)\}$ and $\{s'_1,\ldots,s'_{r'}\}$ are free generating sets of $\mathcal F_K(r')$. Cohn's invariance result \cite[Corollary~5.8.14, p.~258]{C} implies $r=r'$. Finally, since $r\ge n\ge2$, it follows from Cohn's theorem on the center of a free field \cite[Theorem~5.5.11, p.~238]{C} that
$$
Z(D_K(E,\omega))=\Phi\bigl(Z(\mathcal F_K(r))\bigr)=\Phi(K)=K.
$$
\end{proof}

\begin{proposition}
\label{prop:attainable-ranks}
For fixed $n\ge2$, the numbers of free generators arising from LV-roses with $n$ loops whose algebras embed in a division ring are exactly
$$
\{r\in\mathbb N:3n-2\le r\le n^2,\quad r\equiv n\pmod2\}.
$$
As $n$ varies, the possible values of $r(E,\omega)$ are exactly
$$
\{4,7,9,10\}\cup\{r\in\mathbb N:r\ge12\}.
$$
\end{proposition}

\begin{proof}
For such a rose, it follows from the embedding criterion in Theorem~\ref{thm:embedding} that $w_n=n$ and $w_j>j$ for $j<n$. Since the weights are ordered integers, $j+1\le w_j\le n$ for $j<n$. Therefore
$$
\sum_{j=1}^{n-1}(j+1)+n
=\frac{n(n+1)}2+n-1
\le\sum_{j=1}^n w_j\le n^2.
$$
Substituting in $r(E,\omega)=2\sum_jw_j-n^2$, we obtain
$$
3n-2\le r(E,\omega)\le n^2
\quad\text{and}\quad
r(E,\omega)\equiv n\pmod2.
$$

To realize every value in this range with the stated parity, start with the tuple defined by $w_i^{(0)}=i+1$ for $i<n$ and $w_n^{(0)}=n$. If $w\ne(n,\ldots,n)$, choose $j=\max\{i:w_i<n\}$ and replace $w$ by $w+\mathsf{e}_j$, where $\mathsf{e}_j$ is the $j$th standard basis vector of $\mathbb Z^n$. Since $w_{n-1}=w_n=n$, we have $j\le n-2$. The only modified coordinate satisfies
$$
j+1\le w_j<w_j+1\le n=w_{j+1},
$$
and, if $j>1$, also $w_{j-1}\le w_j<w_j+1$. Thus the weights remain nondecreasing, $w_i>i$ for $i<n$, and the last two weights remain $n$. The number of increments from $w^{(0)}$ to $(n,\ldots,n)$ is
$$
h_n=n^2-\left(\frac{n(n+1)}2+n-1\right)
=\frac{(n-1)(n-2)}2.
$$
Writing $w^{(t)}$ for the successive tuples, we have
$$
2\sum_{i=1}^n w_i^{(t)}-n^2=3n-2+2t\quad\text{for }0\le t\le h_n
\quad\text{and}\quad 3n-2+2h_n=n^2.
$$
This proves the assertion for fixed $n$.

Put $\mathcal R_n=\{3n-2+2t:t\in\N,\ 0\le t\le h_n\}$. For $n\ge4$,
$$
\min \mathcal R_{n+2}=3(n+2)-2\le n^2+2=\max \mathcal R_n+2,
$$
and $\mathcal R_n,\mathcal R_{n+2}$ consist of integers of the same parity. Since their upper endpoints tend to infinity,
$$
\mathop{\textstyle\bigcup}\limits_{\substack{n\ge4\\n\text{ even}}}\mathcal R_n
=\{2a:a\in\mathbb N,\ a\ge5\}
\quad\text{and}\quad
\mathop{\textstyle\bigcup}\limits_{\substack{n\ge5\\n\text{ odd}}}\mathcal R_n
=\{2a+1:a\in\mathbb N,\ a\ge6\}.
$$
Together with $\mathcal R_2=\{4\}$ and $\mathcal R_3=\{7,9\}$, these equalities imply
$$
\mathop{\textstyle\bigcup}\limits_{n\ge2}\mathcal R_n=\{4,7,9,10\}\cup\{r\in\mathbb N:r\ge12\}.
$$
\end{proof}

\begin{samepage}
\begin{example}
\label{ex:same-division-different-algebras}
The LV-algebras of the roses with weights $(4,4,4,4)$ and $(2,3,4,5,6,6)$ have $K$-isomorphic universal division rings, both isomorphic to $\mathcal F_K(16)$, since
$$
2(4+4+4+4)-4^2=16
\quad\text{and}\quad
2(2+3+4+5+6+6)-6^2=16.
$$
\begin{center}
\begin{tikzpicture}[>={Stealth[length=1.3mm,width=.9mm]}]
\begin{scope}[shift={(0,0)}]
    \node[inner sep=1pt] at (0,0) (vfour) {$\bullet$};
    \node at (0,-0.55) {$v$};
    \draw[->,rotate=135] (vfour) .. controls (0.30,-0.12) and (1.65,-1.10) .. (1.65,0) .. controls (1.65,1.10) and (0.30,0.12) .. (vfour);
    \node[above left] at (135:1.80) {$a_1,\ldots,a_4$};
    \draw[->,rotate=45] (vfour) .. controls (0.30,-0.12) and (1.65,-1.10) .. (1.65,0) .. controls (1.65,1.10) and (0.30,0.12) .. (vfour);
    \node[above right] at (45:1.80) {$b_1,\ldots,b_4$};
    \draw[->,rotate=225] (vfour) .. controls (0.30,-0.12) and (1.65,-1.10) .. (1.65,0) .. controls (1.65,1.10) and (0.30,0.12) .. (vfour);
    \node[below left] at (225:1.80) {$c_1,\ldots,c_4$};
    \draw[->,rotate=315] (vfour) .. controls (0.30,-0.12) and (1.65,-1.10) .. (1.65,0) .. controls (1.65,1.10) and (0.30,0.12) .. (vfour);
    \node[below right] at (315:1.80) {$d_1,\ldots,d_4$};
    \node at (0,-2.55) {$(4,4,4,4)$};
\end{scope}
\begin{scope}[shift={(6.8,0)}]
    \node[inner sep=1pt] at (0,0) (vsix) {$\bullet$};
    \draw[->,rotate=90] (vsix) .. controls (0.30,-0.12) and (1.65,-0.80) .. (1.65,0) .. controls (1.65,0.80) and (0.30,0.12) .. (vsix);
    \node[above] at (90:1.80) {$e_1,e_2$};
    \draw[->,rotate=150] (vsix) .. controls (0.30,-0.12) and (1.65,-0.80) .. (1.65,0) .. controls (1.65,0.80) and (0.30,0.12) .. (vsix);
    \node[left] at (150:1.80) {$f_1,\ldots,f_3$};
    \draw[->,rotate=210] (vsix) .. controls (0.30,-0.12) and (1.65,-0.80) .. (1.65,0) .. controls (1.65,0.80) and (0.30,0.12) .. (vsix);
    \node[left] at (210:1.80) {$g_1,\ldots,g_4$};
    \draw[->,rotate=270] (vsix) .. controls (0.30,-0.12) and (1.65,-0.80) .. (1.65,0) .. controls (1.65,0.80) and (0.30,0.12) .. (vsix);
    \node[below] at (270:1.80) {$h_1,\ldots,h_5$};
    \draw[->,rotate=330] (vsix) .. controls (0.30,-0.12) and (1.65,-0.80) .. (1.65,0) .. controls (1.65,0.80) and (0.30,0.12) .. (vsix);
    \node[right] at (330:1.80) {$k_1,\ldots,k_6$};
    \draw[->,rotate=30] (vsix) .. controls (0.30,-0.12) and (1.65,-0.80) .. (1.65,0) .. controls (1.65,0.80) and (0.30,0.12) .. (vsix);
    \node[right] at (30:1.80) {$\ell_1,\ldots,\ell_6$};
    \node at (0,-2.55) {$(2,3,4,5,6,6)$};
\end{scope}
\end{tikzpicture}
\end{center}
Only the first rose has all weights equal to its number of loops. The semifir criterion in Corollary~\ref{cor:semifir} therefore implies that exactly one of the two algebras is a semifir, so they are nonisomorphic.
\end{example}
\end{samepage}

\begin{samepage}
\begin{example}
\label{ex:same-monoid-different-division}
The monoid $V(A)$ together with its distinguished class $[A]$ does not determine the universal division ring.
\begin{center}
\begin{tikzpicture}[>={Stealth[length=1.3mm,width=.9mm]}]
\begin{scope}[shift={(0,0)}]
    \node[inner sep=1pt] at (0,0) (vthree) {$\bullet$};
    \node at (0,-0.40) {$v$};
    \draw[->,rotate=90] (vthree) .. controls (0.30,-0.12) and (1.65,-1.20) .. (1.65,0) .. controls (1.65,1.20) and (0.30,0.12) .. (vthree);
    \node[above] at (90:1.80) {$e_1,e_2$};
    \draw[->,rotate=210] (vthree) .. controls (0.30,-0.12) and (1.65,-1.20) .. (1.65,0) .. controls (1.65,1.20) and (0.30,0.12) .. (vthree);
    \node[below left] at (210:1.80) {$f_1,\ldots,f_3$};
    \draw[->,rotate=330] (vthree) .. controls (0.30,-0.12) and (1.65,-1.20) .. (1.65,0) .. controls (1.65,1.20) and (0.30,0.12) .. (vthree);
    \node[below right] at (330:1.80) {$g_1,\ldots,g_3$};
    \node at (0,-2.55) {$A:\ (2,3,3)$};
\end{scope}
\begin{scope}[shift={(6.8,0)}]
    \node[inner sep=1pt] at (0,0) (vfourmixed) {$\bullet$};
    \node at (0,-0.55) {$v$};
    \draw[->,rotate=135] (vfourmixed) .. controls (0.30,-0.12) and (1.65,-1.10) .. (1.65,0) .. controls (1.65,1.10) and (0.30,0.12) .. (vfourmixed);
    \node[above left] at (135:1.80) {$a_1,a_2$};
    \draw[->,rotate=45] (vfourmixed) .. controls (0.30,-0.12) and (1.65,-1.10) .. (1.65,0) .. controls (1.65,1.10) and (0.30,0.12) .. (vfourmixed);
    \node[above right] at (45:1.80) {$b_1,\ldots,b_4$};
    \draw[->,rotate=225] (vfourmixed) .. controls (0.30,-0.12) and (1.65,-1.10) .. (1.65,0) .. controls (1.65,1.10) and (0.30,0.12) .. (vfourmixed);
    \node[below left] at (225:1.80) {$c_1,\ldots,c_4$};
    \draw[->,rotate=315] (vfourmixed) .. controls (0.30,-0.12) and (1.65,-1.10) .. (1.65,0) .. controls (1.65,1.10) and (0.30,0.12) .. (vfourmixed);
    \node[below right] at (315:1.80) {$d_1,\ldots,d_4$};
    \node at (0,-2.55) {$A^{\prime}:\ (2,4,4,4)$};
\end{scope}
\end{tikzpicture}
\end{center}
For the rose with weights $(2,3,3)$, the presentation in Proposition~\ref{prop:monoid} is
$$
V(A)\cong\langle v,q\mid2v=v+q\rangle.
$$
For the rose with weights $(2,4,4,4)$, the corresponding presentation is
$$
V(A')\cong\langle v,q\mid2v=q+v,\ q+2v=3v\rangle.
$$
The second relation follows from the first by adding $v$, so there is a monoid isomorphism $\varphi\colon V(A)\to V(A')$ satisfying $\varphi([A])=[A']$. The universal division rings of $A$ and $A'$ are $\mathcal F_K(7)$ and $\mathcal F_K(12)$, respectively. Since $7\ne12$, Corollary~\ref{cor:field-classification} implies that these division rings are not $K$-isomorphic.
\end{example}
\end{samepage}

\section{Matrices and rational expressions}
\label{sec:matrices}

Let $M\in M_m(A)$. If $M$ factors through a finitely generated projective $A$-module $H$, Corollary~\ref{cor:stablefree} implies that $H\oplus A^s\cong A^{\rho(H)+s}$ for some $s\ge0$. Thus $M\oplus I_s$ factors through $A^{\rho(H)+s}$. Following Cohn \cite[\S0.1, p.~5]{Cohn2006}, we call $M$ \emph{stably full} if $M\oplus I_s$ is full in the sense of Section~\ref{sec:criterion} for every $s\ge0$.

\begin{proposition}
\label{prop:stablefull}
For every $M\in M_m(A)$, the following conditions are equivalent:
\begin{enumerate}
\renewcommand{\theenumi}{\roman{enumi}}
\renewcommand{\labelenumi}{\textup{(\theenumi)}}
\item $r_\rho(M)=m$.
\item $M$ is stably full.
\item $M$ is invertible over $D_K(E,\omega)$.
\end{enumerate}
\end{proposition}

\begin{proof}
Suppose that $M=\gamma\beta$ factors through a finitely generated projective module $H$, with $h=\rho(H)<m$. For some $t\ge0$, there is an isomorphism $\sigma\colon H\oplus A^t\to A^{h+t}$ as in Corollary~\ref{cor:stablefree}. Let $B\colon A^{m+t}\to A^{h+t}$ and $C\colon A^{h+t}\to A^{m+t}$ be the maps defined by
$$
B=\sigma(\beta\oplus\id_{A^t})
\quad\text{and}\quad
C=(\gamma\oplus\id_{A^t})\sigma^{-1}.
$$
Then $CB=(\gamma\oplus\id_{A^t})(\beta\oplus\id_{A^t})=M\oplus I_t$, with $h+t<m+t$. Thus $M$ is not stably full.

Conversely, suppose that for some $s\ge0$ and $q<m+s$ there is a factorization
$$
\begin{pmatrix}M&0\\0&I_s\end{pmatrix}
=\begin{pmatrix}C_1\\C_2\end{pmatrix}\begin{pmatrix}B_1&B_2\end{pmatrix}
$$
through $A^q$, where $B_1\colon A^m\to A^q$, $B_2\colon A^s\to A^q$, $C_1\colon A^q\to A^m$ and $C_2\colon A^q\to A^s$. Comparing blocks, we obtain
$$
C_1B_1=M,\qquad C_2B_1=0\quad\text{and}\quad C_2B_2=I_s.
$$
Put $H=\ker C_2$. Let $\phi\colon H\oplus A^s\to A^q$ and $\psi\colon A^q\to H\oplus A^s$ be the maps defined by
$$
\phi(h,z)=h+B_2z
\quad\text{and}\quad
\psi(a)=(a-B_2C_2a,C_2a).
$$
The first component of $\psi(a)$ lies in $H$, since $C_2(a-B_2C_2a)=C_2a-C_2a=0$. Moreover,
$$
\begin{aligned}
\psi\phi(h,z)&=(h+B_2z-B_2C_2(h+B_2z),C_2h+C_2B_2z)=(h,z),\\
\phi\psi(a)&=a-B_2C_2a+B_2C_2a=a.
\end{aligned}
$$
Thus $H\oplus A^s\cong A^q$, so $H$ is finitely generated projective and $\rho(H)=q-s<m$. Since $C_2B_1=0$, let $\beta_0\colon A^m\to H$ and $\gamma_0\colon H\to A^m$ be the maps defined by $\beta_0(x)=B_1x$ and $\gamma_0=C_1|_H$. Then
$$
A^m\xrightarrow{\ \beta_0\ }H\xrightarrow{\ \gamma_0\ }A^m,
\quad\text{with}\quad \gamma_0\beta_0(x)=C_1B_1x=Mx.
$$ Therefore $r_\rho(M)\le\rho(H)<m$. This proves the equivalence of \textup{(i)} and \textup{(ii)}. Using the rank formula in Lemma~\ref{lem:cohn} for the second equivalence, we obtain
$$
M\in\operatorname{GL}_m(D_K(E,\omega))\quad\Longleftrightarrow\quad
\rk_{D_K(E,\omega)}(M)=m\quad\Longleftrightarrow\quad r_\rho(M)=m,
$$
which proves the remaining equivalence.
\end{proof}

\begin{proposition}
\label{prop:bounded-stabilization}
Put
$$
\kappa=\max_{1\le i<k}\frac{N_i}{d_i-N_i}
\quad\text{and}\quad b_m=\lfloor(m-1)\kappa\rfloor,
$$
with $\kappa=0$ if $k=1$. For every $M\in M_m(A)$, with $m\ge1$,
$$
M\in\Sigma_\rho
\quad\Longleftrightarrow\quad
M\oplus I_{b_m}\text{ is full}.
$$
\end{proposition}

\begin{proof}
Suppose that $M\in\Sigma_\rho$. For every $b\ge0$, $M\oplus I_b$ is invertible over $D_K(E,\omega)$. If $M\oplus I_b=CB$ factors through $A^q$, with $q<m+b$, extension to $D_K(E,\omega)$ would imply
$$
m+b=\rk_{D_K(E,\omega)}(M\oplus I_b)=\rk_{D_K(E,\omega)}(CB)\le q<m+b.
$$
This contradiction proves that $M\oplus I_{b_m}$ is full.

Conversely, suppose that $M\notin\Sigma_\rho$. Choose $M=\gamma\beta$ with $\beta\colon A^m\to H$, $\gamma\colon H\to A^m$, and $h=\rho(H)\le m-1$. Write $[H]=cv+\sum_{i<k}c_iq_i$ and put $s=\sum_{i<k}c_iN_i$. Then
$$
h=c+\sum_{i<k}c_i(d_i-N_i)
\quad\text{and}\quad
s\le\kappa\sum_{i<k}c_i(d_i-N_i)\le\kappa h\le\kappa(m-1).
$$
Since $s$ is an integer, $s\le b_m$. For this value of $s$, Corollary~\ref{cor:stablefree} ensures that there is an isomorphism $\sigma\colon H\oplus A^s\to A^{h+s}$. Put $b=b_m$, and let $\sigma_b\colon H\oplus A^b\to A^{h+b}$ be the isomorphism defined by $\sigma_b=\sigma\oplus\id_{A^{b-s}}$.
Let $B\colon A^{m+b}\to A^{h+b}$ and $C\colon A^{h+b}\to A^{m+b}$ be the maps defined by
$$
B=\sigma_b(\beta\oplus\id_{A^b})
\quad\text{and}\quad
C=(\gamma\oplus\id_{A^b})\sigma_b^{-1}.
$$
Then
$$
CB=(\gamma\beta)\oplus I_b=M\oplus I_b,
\quad\text{with}\quad h+b<m+b.
$$
Thus $M\oplus I_{b_m}$ is not full.
\end{proof}

Following Cohn \cite[\S6.6, pp.~317--318]{C}, we say that a presentation has a \emph{solvable word problem} if an algorithm decides whether any two well-defined expressions in its generators represent the same element.

For a presented division ring, the task is to find an effective test for $u-v=0$, where $u$ and $v$ are rational expressions defined in that division ring. For a semifir $R$, Cohn reduces the word problem in its universal field of fractions to deciding whether square matrices over $R$ are full \cite[Theorem~6.6.2, p.~319]{C}. For a free field over a commutative field, he proves that the word problem is solvable whenever the coefficient field is given by a presentation with solvable word problem \cite[Corollary~6.6.8, p.~323]{C}.

The isomorphism constructed in Theorem~\ref{thm:free-coordinates} identifies $D_K(E,\omega)$ with a free field and expresses the entries of $X$ and $Y$ as rational expressions in the free generators. If the given presentation of $K$ has solvable word problem, these formulas allow us to apply Cohn's zero test to rational expressions over $A$. The same test selects pivots in Gaussian elimination over the free field. The rank formula in Lemma~\ref{lem:cohn} identifies the rank computed for $M$ with $r_\rho(M)$.

\begin{corollary}
\label{cor:effective}
Suppose that $K$ is given by a presentation with solvable word problem. Then one can decide whether two rational expressions in $D_K(E,\omega)$ are equal and whether a square matrix over $A$ is $\rho$-full.
\end{corollary}

\begin{proof}
Put $r=r(E,\omega)$ and $F=\mathcal F_K(r)$. Let $\Phi\colon F\to D_K(E,\omega)$ be the isomorphism of Theorem~\ref{thm:free-coordinates}, and let $\tau\colon D_K(E,\omega)\to F$ be the isomorphism defined by $\tau=\Phi^{-1}$. The free generators of $F$ are the images under $\tau$ of the elements $p_j$ for $1\le j\le n$ and $l_{ij},t_{ji}$ for $(i,j)\in\mathcal I$, and
$$
\tau(X)=\tau(L)\tau(P)\tau(T)^{-1}
\quad\text{and}\quad
\tau(Y)=\tau(T)\tau(P)^{-1}\tau(L)^{-1}.
$$
Here the unitriangular inverses are finite polynomial expressions, and $\tau(P)^{-1}$ is diagonal with entries the inverses of the corresponding free generators. These formulas compute $\tau$ on the generators of $A$. For a rational expression over $A$, apply recursively
$$
\tau(a+b)=\tau(a)+\tau(b),\qquad
\tau(ab)=\tau(a)\tau(b)\quad\text{and}\quad
\tau(a^{-1})=\tau(a)^{-1}\quad\text{whenever }a\ne0.
$$
The given presentation of $K$ has solvable word problem, so arithmetic in $K$ is computable. It follows from Cohn's result for free fields \cite[Corollary~6.6.8, p.~323]{C} that one can decide whether $\tau(x)-\tau(y)=0$, equivalently $x=y$ in $D_K(E,\omega)$.

For $M\in M_m(A)$, Gaussian elimination over $F$ uses these zero tests to select pivots and computes $U,V\in\operatorname{GL}_m(F)$ such that
$$
U\tau(M)V=\operatorname{diag}(I_q,0)\qquad\text{with }0\le q\le m.
$$
Consequently $\rk_{D_K(E,\omega)}(M)=\rk_F(\tau(M))=q$. In view of the rank identity in Lemma~\ref{lem:cohn},
$$
M\in\Sigma_\rho\quad\Longleftrightarrow\quad r_\rho(M)=m
\quad\Longleftrightarrow\quad q=m.
$$
\end{proof}

\begin{proposition}
\label{prop:fractions}
Every $x\in D_K(E,\omega)$ has an expression $x=uM^{-1}v$, where $M\in\Sigma_\rho\cap M_m(A)$, $u\in M_{1\times m}(A)$ and $v\in M_{m\times1}(A)$. If $y=pN^{-1}q$ is another such expression, then
$$
x+y=\begin{pmatrix}u&p\end{pmatrix}
\begin{pmatrix}M&0\\0&N\end{pmatrix}^{-1}
\begin{pmatrix}v\\q\end{pmatrix}
\quad\text{and}\quad
xy=\begin{pmatrix}u&0\end{pmatrix}
\begin{pmatrix}M&-vp\\0&N\end{pmatrix}^{-1}
\begin{pmatrix}0\\q\end{pmatrix}.
$$
For $B=\begin{pmatrix}0&u\\v&M\end{pmatrix}$, we have $x=0$ if and only if $B\notin\Sigma_\rho$. If $x\ne0$, then $x^{-1}=-\eta_1^TB^{-1}\eta_1$, where $\eta_1=(1,0,\ldots,0)^T\in A^{m+1}$.
\end{proposition}

\begin{proof}
The two block matrices in the sum and product have inverses
$$
\begin{pmatrix}M&0\\0&N\end{pmatrix}^{-1}
=\begin{pmatrix}M^{-1}&0\\0&N^{-1}\end{pmatrix}
\quad\text{and}\quad
\begin{pmatrix}M&-vp\\0&N\end{pmatrix}^{-1}
=\begin{pmatrix}M^{-1}&M^{-1}vpN^{-1}\\0&N^{-1}\end{pmatrix}.
$$
Since both block matrices are invertible over $D_K(E,\omega)$, they belong to $\Sigma_\rho$ by the criterion in Proposition~\ref{prop:stablefull}. The indicated products are
$$
\begin{pmatrix}u&p\end{pmatrix}\begin{pmatrix}M^{-1}&0\\0&N^{-1}\end{pmatrix}\begin{pmatrix}v\\q\end{pmatrix}
=uM^{-1}v+pN^{-1}q=x+y
$$
and
$$
\begin{pmatrix}u&0\end{pmatrix}\begin{pmatrix}M^{-1}&M^{-1}vpN^{-1}\\0&N^{-1}\end{pmatrix}\begin{pmatrix}0\\q\end{pmatrix}
=uM^{-1}vpN^{-1}q=xy.
$$
Let $\mathcal R\subseteq D_K(E,\omega)$ consist of the elements $uM^{-1}v$ with $M,u,v$ as in the statement. The formulas imply $\mathcal R+\mathcal R\subseteq\mathcal R$ and $\mathcal R\mathcal R\subseteq\mathcal R$. In addition,
$$
a=a(1)^{-1}1,\qquad
(M^{-1})_{ij}=\mathsf{e}_i^TM^{-1}\mathsf{e}_j
\quad\text{and}\quad
-(uM^{-1}v)=(-u)M^{-1}v,
$$
where $\mathsf{e}_i=(\delta_{1i},\ldots,\delta_{mi})^T\in A^m$ for $1\le i\le m$. Hence $\mathcal R$ is a subring containing $A$ and the entries of $M^{-1}$ for all $M\in\Sigma_\rho$. By the presentation in Section~\ref{sec:construction}, these entries and $A$ generate $D_K(E,\omega)$ as a ring. Thus $\mathcal R=D_K(E,\omega)$.

For $x=uM^{-1}v$, put $B=\begin{pmatrix}0&u\\v&M\end{pmatrix}$. Direct multiplication verifies
$$
\begin{pmatrix}1&-uM^{-1}\\0&I_m\end{pmatrix}
B
\begin{pmatrix}1&0\\-M^{-1}v&I_m\end{pmatrix}
=\begin{pmatrix}-x&0\\0&M\end{pmatrix}.
$$
Therefore
$$
B\in\Sigma_\rho\quad\Longleftrightarrow\quad
B\in\operatorname{GL}_{m+1}(D_K(E,\omega))\quad\Longleftrightarrow\quad x\ne0.
$$
For $x\ne0$, the inverse identity is
$$
B^{-1}=\begin{pmatrix}1&0\\-M^{-1}v&I_m\end{pmatrix}
\begin{pmatrix}-x^{-1}&0\\0&M^{-1}\end{pmatrix}
\begin{pmatrix}1&-uM^{-1}\\0&I_m\end{pmatrix}.
$$
In particular, $\eta_1^TB^{-1}\eta_1=-x^{-1}$, which proves the inverse formula.
\end{proof}

The representations are not unique, since
$$
uM^{-1}v=\begin{pmatrix}u&0\end{pmatrix}
(M\oplus I_1)^{-1}\begin{pmatrix}v\\0\end{pmatrix}.
$$
To test whether $x=y$, write $x-y=uM^{-1}v$ and put $B=\begin{pmatrix}0&u\\v&M\end{pmatrix}$. The zero criterion in Proposition~\ref{prop:fractions} says that $x=y$ if and only if $B\notin\Sigma_\rho$. If $K$ is given by a presentation with solvable word problem, the $\rho$-fullness test in Corollary~\ref{cor:effective} decides whether $B$ belongs to $\Sigma_\rho$.

\begin{example}
\label{ex:mixed}
Let $(E,\omega)$ be the LV-rose with weights $(2,3,3)$ in Example~\ref{ex:same-monoid-different-division}, and put $A=L_K(E,\omega)$. Then
$$
V(A)\cong\langle v,q\mid2v=v+q\rangle
\quad\text{with}\quad\rho(v)=\rho(q)=1.
$$
Since $w_3=3$ and $w_j>j$ for $j<3$, the embedding criterion in Theorem~\ref{thm:embedding} implies that $A$ embeds in a division ring. As in the proof of Corollary~\ref{cor:semifir}, $A$ cannot be a semifir because $q\ne mv$ for every $m\in\N$. Let $e,f,g$ be its loops, with $e$ of weight $2$. Then
$$
X=\begin{pmatrix}e_1&f_1&g_1\\e_2&f_2&g_2\\0&f_3&g_3\end{pmatrix}
\quad\text{and}\quad
Y=\begin{pmatrix}e_1^*&e_2^*&0\\f_1^*&f_2^*&f_3^*\\g_1^*&g_2^*&g_3^*\end{pmatrix}.
$$
Put $u=(e_1,e_2)^T$, $u^*=(e_1^*,e_2^*)$, and $\epsilon=I_2-uu^*$. We have $u^*u=1$, $\epsilon^2=\epsilon$, and $P=\epsilon A^2$ satisfies
$$
[P]=q,\qquad P\oplus A\cong A^2\quad\text{and}\quad P\not\cong A.
$$
The matrix $\epsilon$ is full. Indeed, it is nonzero since $P\ne0$. If $\epsilon=ab$ with $a\in M_{2\times1}(A)$ and $b\in M_{1\times2}(A)$, then
$$
a(ba-1)b=\epsilon^2-\epsilon=0.
$$
The algebra $A$ is a domain because it embeds in $D_K(E,\omega)$. For nonzero entries $a_i$ and $b_j$, the $(i,j)$-entry of the last equality is $a_i(ba-1)b_j=0$. Hence $ba=1$. The maps from $A$ to $P$ and from $P$ to $A$ defined by $z\mapsto az$ and $p\mapsto bp$, respectively, would then be inverse, contrary to $P\not\cong A$. Therefore $\epsilon$ has no factorization through $A^1$.

In contrast, $\epsilon u=0$, $u^*\epsilon=0$ and $u^*u=1$ imply
$$
\epsilon\oplus I_1=\begin{pmatrix}\epsilon\\u^*\end{pmatrix}
\begin{pmatrix}\epsilon&u\end{pmatrix}.
$$
This factors a $3\times3$ matrix through $A^2$, so $\epsilon\oplus I_1$ is not full. Thus $\epsilon$ is not stably full, and $r_\rho(\epsilon)=\rho(P)=1<2$. Thus $\epsilon$ is full, but $\epsilon\notin\Sigma_\rho$.
\end{example}

\section*{Funding}
H.~V.~Khanh is supported by the Vietnam National Foundation for Science and Technology Development (NAFOSTED) under Grant No.~101.04-2025.41.

\end{document}